\documentclass[10pt]{elsarticle}

\usepackage[english]{babel}
\usepackage{amsfonts}
\usepackage{amsmath}
\usepackage{amssymb}
\usepackage{graphicx}
\usepackage{color,bm}
\usepackage{latexsym}
\usepackage{epstopdf}
\usepackage{subfigure}
\usepackage[colorlinks=true,linkcolor=blue,citecolor=red,urlcolor=cyan]{hyperref}

\usepackage[normalem]{ulem}

\journal{.}

\newcommand{\bnabla}	{{\mathbf{\nabla}}}
\newcommand{\ddroit}	{{\mathrm{d}}}
\newcommand{\bzero}		{{\mathbf{0}}}
\newcommand{\transpose}	{\mathsf{T}}
\newcommand{\bstate}    {\mathbf{v}}
\newcommand{\nstate}    {N}

\newcommand{\bsens}      {\mathbf{y}}
\newcommand{\nobs}      {n}
\newcommand{\bpolicy}   {\mathbf{\pi}}
\newcommand{\bomega}    {\boldsymbol{\omega}}
\newcommand{\btheta}    {\boldsymbol{\theta}}

\begin{document}

\author{Michele Alessandro Bucci\fnref{label1}}
\ead{michelealessandro.bucci [at] limsi.fr}
\author{Onofrio Semeraro\fnref{label1}}
\author{Alexandre Allauzen\fnref{label1}}
\author{Guillaume Wisniewski\fnref{label1}}
\author{Laurent Cordier\fnref{label2}}
\author{Lionel Mathelin\fnref{label1}}

\address[label1]{LIMSI, CNRS, Universit\'e de Paris-Saclay, Orsay, FR}
\address[label2]{Institut Pprime, CNRS, Universit\'e de Poitiers, ISAE-ENSMA, Poitiers, FR}

\title{Control of chaotic systems by Deep Reinforcement Learning}

%============================================================
% body of the paper
%============================================================

\begin{abstract}
Deep Reinforcement Learning (DRL) is applied to control a nonlinear, chaotic system governed by the one-dimensional Kuramoto-Sivashinsky (KS) equation.
DRL uses reinforcement learning principles for the determination of optimal control solutions and deep Neural Networks for approximating the value function and the control policy. Recent applications have shown that DRL may achieve superhuman performance in complex cognitive tasks.

\smallskip
In this work, we show that using restricted, localized actuations, partial knowledge of the state based on limited sensor measurements, and model-free DRL controllers, it is possible to stabilize the dynamics of the KS system around its unstable fixed solutions, here considered as target states. The robustness of the controllers is tested by considering several trajectories in the phase-space emanating from different initial conditions; we show that the DRL is always capable of driving and stabilizing the dynamics around the target states. 

\smallskip
The complexity of the KS system, the possibility of defining the DRL control policies by solely relying on the local measurements of the system, and their efficiency in controlling its nonlinear dynamics pave the way for the application of RL methods in control of complex fluid systems such as turbulent boundary layers, turbulent mixers or multiphase flows.
\end{abstract}

\maketitle

\section{Introduction}  \label{Sec_Intro}

With the availability of unprecedented computational resources,  the maturity of modeling in many areas of engineering has yielded systems close to optimality in the context of well-known settings they were engineered for.  To further improve such systems, a scientific challenge lies in rendering these systems adaptive by design to changes in the operating conditions. Enlarging the perspective from the engineering standpoint, the current environmental needs have also invigorated the research effort on \emph{flow control} applications. For example, carbon dioxide emissions are considered one of the main responsible for the global warming and any reduction of these emissions can lead to an attenuation of this effect. Increasing the efficiency of the existing technologies for the production of sustainable energy can lead to high potential benefits or to a substantial reduction of oil-consumption in the transport economy sector. Not surprisingly, a rather large body of literature is already available on the many different efforts in this realm, with varying assumptions on the system, such as linearity of the governing equations, or on the amount of information one has on the system, such as observability of a state vector or not. The interested reader is referred to the reviews by \cite{bewley2007linear,sipp2016amr,Brunton_Noack_review}.

Active control for the optimization of the performance can be introduced via adequate strategies capable of modifying the system response in a prescribed manner (\emph{open-loop}) or as a function of some observations of the system at hand (\emph{closed-loop}). In both cases, the challenge is to infer an efficient and robust control strategy, hereafter termed \emph{policy}, improving upon the retained objective function. In the usual \emph{model-based} approach, a dynamical model is used to describe the behavior of the system.

This model allows to predict the effect of a given control action and can hence be used to derive the best control strategy leading to an optimal performance.  However, a physical model\footnote{The term \emph{model} is ambiguous. Here, the term \textit{model} refers to a \emph{physical model}. This meaning differs from its use in the Machine Learning community where the term is more related to a parameterized function that links inputs to outputs.} is not always available. Besides systems for which the governing equations are simply unknown or very poorly known, there are many situations where solving the governing equations is too slow with respect to the dynamics at play to be useful. While reduced-order models may help in solving an approximate system meeting real-time constraints, they usually lack robustness and can critically lose accuracy when control is applied, resulting in poor performance, at best.

A different line of control strategy relies on a \emph{data-driven} approach. In this view, no model is employed and the control command is derived based on past observations only. In a training step, control actions are applied to the system and observations are made, possibly including the evaluation of the objective function. From this collection of observations, an input-output (I/O) model is built and subsequently used for controlling the system. Typical of this viewpoint are extremum seeking, \cite{Becker_etal_07}, control strategies relying on system-identification techniques leading to auto-regressive models (AR, ARMA(X), etc.), \textit{e.g.}, \cite{Kegerise_01, Kegerise_02}, or the more recently proposed so-called Machine Learning Control, \cite{gautier_aider_duriez_noack_segond_abel_2015}. These techniques are more directly compatible with real systems since they do not require a physical model, nor prohibitive computational power. Moreover, they do not assume knowledge of a state vector and can efficiently deal with partial observations. Yet, they either exhibit limited performance in the general case where prior expertise knowledge is not available, or require an extensive training set of trial-and-errors, negatively affecting the learning time.
\begin{figure}
\begin{center}
    \includegraphics[width=10.0cm]{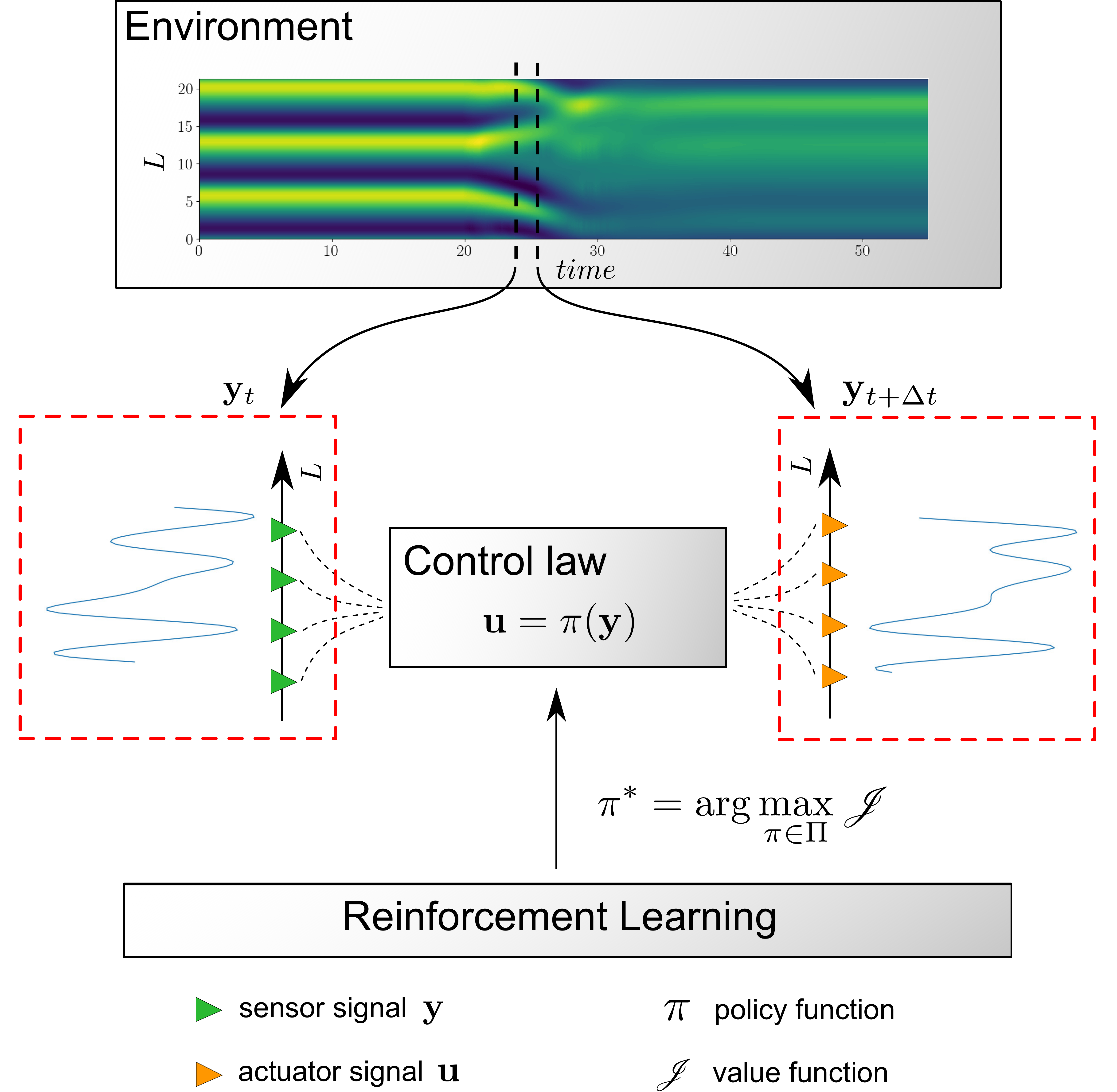}
    \caption{    
Overview of the Reinforcement Learning (RL) strategy for the Kuramoto-Sivashinsky equation considered in Sec.~\ref{sec:plant}.
In the RL framework, an \emph{agent} interacts with an \emph{environment} by making \emph{observations} $\bsens$ and performing \emph{actions} $\mathbf{u}$.
In return, the agent receives a \emph{reward} that depends directly on the changes of the environment induced by the action.
The objective of RL is to determine the action $\mathbf{u}$ to impose on the environment in order to maximize a given \emph{value function} $\mathcal{J}$ that represents the cumulative rewards over time. 
The \emph{control law} $\mathbf{u}$ is determined via the \emph{policy function} $\bpolicy$ which maps measurements $\bsens$ taken from the environment to the action space.
}
    \label{fig:bigpicture}
\end{center}
\end{figure}
In the present work, we introduce a Reinforcement Learning (RL) strategy \cite{sutton1998introduction} for the closed-loop, nonlinear control. RL is a well-established technique mainly originating from the robotics community and has gained wide-spread popularity with some recent mediatic applications, achieving super-human performances in the \texttt{go} and \texttt{shogi} games as well as self-teaching for the \texttt{StarCraft II} videogame, \cite{silver2017mastering}, or in revenue management \cite{Gosavi_19} to cite only a few examples. As a data-driven technique, it shares the applicability of other I/O approaches, \textit{e.g.}, low computational requirements, and the ability to use only limited sensors in contrast with a full state vector. Reinforcement learning dates back to the 1950s when the problem of optimal control was solved by Richard Bellman through the introduction of dynamic programming \cite{bellman1958dynamic}, leading -- in the continuous formulation -- to the Hamilton-Jacobi-Bellman (HJB) equation.

Figure~\ref{fig:bigpicture} provides a sketch of the RL methodology for the control of dynamical system. The physical system under consideration, termed the \emph{environment} in the reinforcement learning literature, is observed at time $t$ by a so-called \emph{agent} through a set of localized measurements $\bsens(t) \in \mathbb{R}^{\nobs}$. 
The agent performs an action $\mathbf{u}(t)$ that changes the future state of the environment and, in return, receives a \emph{reward} that represents the degree at which a state/action pair is desirable for the targeted objective function.
This measure of performance to be optimized is then defined as the expectation of the (discounted) cumulative rewards over time.
The action $\mathbf{u}(t) \in \mathbb{R}^m$ is evaluated from the current observations via the \emph{control policy} $\bpolicy$,  defined as $\mathbf{u}(t) = \bpolicy\left(\bsens(t)\right)$, which represents a mapping from the observation space ($\mathbb{R}^{\nobs}$) to the action space ($\mathbb{R}^m$). The policy is defined in a given class $\Pi$ of multivariate functions. 
The \emph{value function} $\mathcal{J}: \mathbb{R}^{\nobs} \rightarrow \mathbb{R}$ or cost-to-go function of a policy $\bpolicy$ over a time horizon gives the cumulative rewards when $\bsens_t$ is the current measurement and the system follows policy $\bpolicy$ thereafter. 
%In DRL, $\mathcal{J}$ is also represented by a neural network (NN). 
The optimal policy, denoted $\bpolicy^\star$, maximizes $\mathcal{J}$ over $\Pi$. 

Estimation of the value function unlocks avenues for a more efficient use of the available information, potentially resulting in a faster training, and improved robustness with respect to small perturbations to the system and to the sensor measurements. Applications of RL have been mostly restricted to the discrete settings, where the number of possible actions is potentially huge, but finite, \cite{Gorodetsky2015_TT}. A notable extension to the continuous framework includes the work by Gorodetsky and collaborators using function trains, \cite{Gorodetsky_FT}. Our earlier efforts in bringing RL to the context of fluid flows involved Q-learning, \cite{Gueniat_etal_16}, and temporal difference-based continuous approaches, \cite{Pivot_etal_17, Pivot_etal_18}. 

Recent efforts from the literature involve optimizing a collective swimming strategy \cite{Verma_etal_18} and the control of the flow around a circular cylinder in the laminar regime \cite{Rabault_JFM_19}, among others. 
Here, we demonstrate the performance of our methodology in terms of quality of the control strategy (achieved performance) and robustness with respect to perturbations to the system and the sensor measurements.

In this proof-of-concept work, we focus on the control of a nonlinear, chaotic dynamical system, the 1D Kuramoto-Sivashinsky (KS) equation as a model of a fluid flow to be controlled. 
This model is described in a space-time domain with a 4-th order partial differential equation (PDE). 
The KS system exhibits some of the typical features observed in flow systems at low Reynolds number, such as traveling waves, and may reach a chaotic regime closer to turbulence for certain configurations, \cite{holmes1996turbulence,bohr2005dynamical,cvitanovic2010state}. Due to these peculiarities, the KS system serves as a challenging test-bed for more complex fluidic systems to be controlled such as turbulent mixers or turbulent boundary layers. 

\bigskip
The paper is organized as follows. The basics of our control strategies are presented in Sec.~\ref{sec:nonopt} and \ref{sec:RL}.
In Sec.~\ref{sec:nonopt}, we review how Reinforcement Learning is connected to the classical optimal control problem.
We first derive the Hamilton-Jacobi-Bellman equation from optimality principles, and next particularize in the discrete time settings to obtain the Bellman equation.
The connection with linear optimal controllers is further described in App.~\ref{sec:appA}. 
In Sec.~\ref{Sec_severalRLs}, we briefly discuss the different classes of RL algorithms.
In Sec.~\ref{sec:ddpg}, the Deep Deterministic Policy Gradient (DDPG) reinforcement learning technique we are using is then presented.
The DRL methodology is illustrated with the control of the Kuramoto-Sivashinsky model in Sec.~\ref{sec:plant} where results are put in perspective with insights from the physics. The paper finalizes with conclusions in Section~\ref{sec:conclusions}. 

% reinforcement learning
\section{From nonlinear optimal control to Reinforcement Learning}\label{sec:nonopt}
The section briefly introduces the mathematical background of Reinforcement Learning. For a deeper introduction to RL, we refer the interested reader to the books by \cite{sutton1998introduction,goodfellow2016deep}. From the methodological viewpoint, we introduce RL following the approach taken in \cite{kirk2012optimal,recht2018tour}, starting from the definition of the optimal control (Sec.~\ref{sec:controdef}) as basis for the Hamilton-Jacobi-Bellman (HJB) equation (Sec.~\ref{sec:HJB}) and stressing the analogies with optimal control theory. Further, we show how we can derive the time-discrete counterpart of the HJB equation, the Bellman equation (Sec.~\ref{sec:bellman}), which serves as the theoretical and mathematical ground for RL applications. 

\subsection{Definition of the control problem}\label{sec:controdef}
We consider a dynamical system to be controlled with localized inputs (\emph{actuators}) and outputs (\emph{sensors}), see for instance Fig.~\ref{fig:bigpicture} where these elements are organized along the streamwise coordinate $x$. Hereafter, the combination of the system to be controlled, the actuators and the sensors will be simply referred to as the \emph{plant}, following the classical terminology in control theory. From the mathematical view point, this corresponds to defining the state-space model that, in the most general case, reads
\begin{subequations}
\begin{align}
\dfrac{\ddroit \bstate}{\ddroit t} &= \mathbf{f}\left(\bstate(t),\mathbf{u}(t),t\right), \label{eq:statespaceNLN1}\\
\bsens(t) &= \mathbf{g}\left(\bstate(t),\mathbf{u}(t),t\right).\label{eq:statespaceNLN2}
\end{align}
\label{eq:statespaceNLN}
\end{subequations}
Equation \eqref{eq:statespaceNLN1} is the state equation, where the map $\mathbf{f}$ propagates in time the state $\bstate\in\mathbb{R}^{\nstate}$, while \eqref{eq:statespaceNLN2} is the output equation. In the optimal control problem \cite{lewis2012optimal}, we aim at defining a control signal $\mathbf{u}\in\mathbb{R}^{m}$ feeding the actuators, based on the sensor measurements $\bsens\in\mathbb{R}^{\nobs}$, such that an objective function $\mathcal{J}$ is minimized 
\begin{equation}
\mathcal{J} = 
h\left(\bstate(T),T\right) +
\int_0^{T} r\left(\bstate(\tau),\mathbf{u}(\tau),\tau\right) \,\ddroit \tau,
\label{eq:objective_function} 
\end{equation}
where $h$ is a specified function, $r$ is a reward associated with the action $\mathbf{u}$ and $T$ is the optimization horizon. According to \cite{kirk2012optimal}, this optimization problem can be embedded in a larger class of problems by considering
\begin{equation}
\mathcal{J}(\bstate_t,t,\underset{t\le\tau\le T}{\mathbf{u}(\tau)}) = 
h\left(\bstate(T),T\right) +
\int_t^{T} r\left(\bstate(\tau),\mathbf{u}(\tau),\tau\right) \,\ddroit \tau,
\label{eq:objective_function_embedded} 
\end{equation}
where $t$ can be any value less than or equal to $T$. Following the convention tacitly introduced in Sec.~\ref{Sec_Intro}, 
we note $\mathcal{J}^\star=\displaystyle\max_{\mathbf{u}}\mathcal{J}$ the optimal value of the objective function. 
The controller (or more properly said the \emph{compensator}) provides the mapping between the measurements $\bsens$ of the system and the control actions $\mathbf{u}$. For now, we stress that the aim of RL is to determine an optimal policy function $\bpolicy^\star$ that describes the optimal control $\mathbf{u}^\star$ from the current observations $\bsens$ \cite{glad2000control, skogestad2005book} following
\begin{eqnarray}
\mathbf{u}^\star(t) = \bpolicy^\star\left(\bsens(t),t\right).
\end{eqnarray}

\subsection{Hamilton-Jacobi-Bellman equation}\label{sec:HJB}
The optimal control problem stated in \eqref{eq:statespaceNLN}-\eqref{eq:objective_function} can be solved by maximizing an augmented Lagrangian \cite{lewis2012optimal} where the governing equations act as constraints. After optimal conditions are imposed on the Lagrangian, a direct-adjoint optimality system can be derived.
When the function $\mathbf{f}$ is linear or can be linearized, and the objective function is quadratic, the final controller is obtained as the solution to a Riccati equation (see Appendix~\ref{sec:appA} for the derivation).
Here, we focus on the general nonlinear case \eqref{eq:objective_function}, and follow a dynamic programming approach to solve the optimal control problem. 
The objective is to determine a Partial Differential Equation (PDE) for the optimal objective function (or value function in the RL terminology) $\mathcal{J}^\star$ such that the corresponding action satisfies the constraint given by the state equation~\eqref{eq:statespaceNLN}.
By definition, the maximum value of the objective function is equal to
\begin{equation}
\label{eq:minimum_cost_function} 
\begin{split}
\mathcal{J}^\star(\bstate(t),t) & = 
\max_{\mathbf{u}}
\mathcal{J}(\bstate_t,t,\underset{t\le\tau\le T}{\mathbf{u}(\tau)})\\
& =
\max_{\underset{t\le\tau\le T}{\mathbf{u}(\tau)}}
\left[
 \int_t^{T} r\left(\bstate(\tau),\mathbf{u}(\tau),\tau\right) \,\ddroit \tau +
h\left(\bstate(T),T\right)
\right].
\end{split}
\end{equation}
By splitting the integrand in \eqref{eq:minimum_cost_function} between the immediate reward in the interval $[t,\, t+\Delta t]$ and the future value function at $t+\Delta t$, we obtain
\begin{equation}
\label{eq:minimum_cost_function2} 
\mathcal{J}^\star(\bstate(t),t) = 
\max_{\underset{t\le\tau\le T}{\mathbf{u}(\tau)}}
\left[
 \int_t^{t+\Delta t} r\,\ddroit \tau +
  \int_{t+\Delta t}^{T} r\,\ddroit \tau +
h\left(\bstate(T),T\right)
\right].
\end{equation}
The principle of optimality requires that
\begin{equation}
\label{eq:minimum_cost_function3} 
\mathcal{J}^\star(\bstate(t),t) = 
\max_{\underset{t\le\tau\le t+\Delta t}{\mathbf{u}(\tau)}}
\left[
 \int_t^{t+\Delta t} r\,\ddroit \tau +
 \mathcal{J}^\star(\bstate(t+\Delta t),t+\Delta t)
\right].
\end{equation}
Expanding $\mathcal{J}^\star(\bstate(t+\Delta t),t+\Delta t)$ in a Taylor series about $(\bstate(t),t)$ and taking the limit as $\Delta t \rightarrow 0$ gives
\begin{equation}
- \dot{\mathcal{J}^\star}(\bstate(t),t) = 
\max_{\mathbf{u}(t)} \left[ r\left(\bstate(t),\mathbf{u}(t),t\right) + \mathcal{J}^\star_{\bstate}\left(\bstate(t),t\right) \, \mathbf{f}\left(\bstate(t),\mathbf{u}(t),t\right) \right],
\label{eq:HJB}
\end{equation}
where $\dot{\mathcal{J}^\star}$ is the temporal derivative of the optimal value function and $\mathcal{J}^\star_{\bstate}$ is the derivative with respect to the state. 
Equation \eqref{eq:HJB} is the Hamilton-Jacobi-Bellman (HJB) equation.
This equation is continuous in time and defined backward. The terminal condition is given by
\begin{equation}
\mathcal{J}^\star (\bstate(T),T) = h\left(\bstate(T),T\right).
\label{eq:HJBTerminal}
\end{equation}

The HJB equation is a sufficient condition for an optimum \cite{bertsekas1996dynamic,Stengel_1994}. Indeed, a value function might fail to satisfy the differentiability and continuity conditions that are required to solve \eqref{eq:HJB} and yet still be optimal.
When solved over the whole state space and when the value function is continuously differentiable, the HJB equation becomes a necessary and sufficient condition for an optimum \cite{kirk2012optimal}.
In case of a continuous action-state space, the optimal policy $\bpolicy^\star$ is the one that produces the optimal trajectory by obeying the HJB. 
However, the whole action-state space is rarely known, especially when the dimension of $\mathbf{f}$ is large.

Note that, for infinite horizon optimizations, it is common to introduce a \emph{discount rate factor} $\rho>0$, that penalizes the immediate reward in the future. In this case, \eqref{eq:minimum_cost_function} becomes
\begin{equation}
\mathcal{J}^\star(\bstate(t),t) = \max_{\mathbf{u}} \int_t^\infty {e^{-\rho \tau} \, r\left(\bstate(\tau),\mathbf{u}(\tau),\tau\right) \, \ddroit \tau}.
\end{equation} 
If we assume that $\mathcal{J}^\star$, $\mathbf{f}$ and $r$ do not explicitly depend on the time $t$, the HJB equation then finally rewrites 
\begin{equation}
\rho \mathcal{J}^\star(\bstate) = \max_{\mathbf{u}} \left[ r(\bstate,\mathbf{u}) + \mathcal{J}^\star_{\bstate}(\bstate) \, \mathbf{f}(\bstate,\mathbf{u}) \right].
\label{eq:HJB_infinite} 
\end{equation}

\subsection{Bellman equation}\label{sec:bellman}
In the derivation of the HJB equation, we have tacitly assumed that the model \eqref{eq:statespaceNLN} is known or can be inferred, for instance by system identification. This is usually not the case as, most of the time, only the state $\bstate$ or a reduced-order representation at a given time $t$ can be accessed, for instance by instantaneous measurements. In this case, the right framework is to consider a Markov Decision Process (MDP) for which the decision making is formulated 
by means of a transition matrix expressing the probability of evolving from a state to another under the chosen action $\mathbf{u} $. 
This framework introduces a discrete time stochastic control process and, as such, requires the reformulation of the objective function in terms of expectation \cite{sutton1998introduction}. Letting $\bstate_t$ and $\mathbf{u}_t$ being the state and the action at the discrete time $t$, respectively, and $\bstate_{t+\Delta t}$ be the state at $t+\Delta t$, \eqref{eq:HJB_infinite} can be rewritten as
\begin{equation}
\mathcal{J}^\star(\bstate_t) = \max_\mathbf{u} \left[ \Delta t \, r(\bstate_t,\mathbf{u}_t) + \gamma \mathcal{J}^\star(\bstate_{t+\Delta t}) \right],
\label{eq:Bellman_opt}
\end{equation}
where $\gamma = \exp\left(-\Delta t \, \rho\right)$ is the discount factor. Including $\Delta t$ in the definition of $r(\bstate_t,\mathbf{u}_t)$, \eqref{eq:Bellman_opt} is the Bellman optimality equation \cite{bellman1958dynamic}.
This equation, which describes the evolution of the optimal value function under the optimal policy $\bpolicy^\star$, is the foundation of the dynamic programming theory.

Let $\mathcal{J}^{\bpolicy}(\bstate_t)$ denote the long-term reward achieved for the state $\bstate_t$, and following a particular policy $\pi$. With the developments made in Sec.~\ref{sec:HJB} for the HJB equation, we derive the Bellman equation for the value function given by
\begin{equation}
\mathcal{J}^{\bpolicy}(\bstate_t) = r(\bstate_t, \mathbf{u}_t) + \gamma \mathcal{J}^{\bpolicy}(\bstate_{t+\Delta t}).
\label{eq:Bellman_value}
\end{equation}

Similarly, we can derive a Bellman equation for the \emph{state-action value function} $Q^{\bpolicy}(\bstate_t, \mathbf{u}_t)$ or \emph{$Q$-function}. This quantity is a measure of the long-term reward assuming the agent is in state $\bstate_t$, performs action $\mathbf{u}_t$, and then continues following some policy $\pi$. The Bellman equation for the Q-function is written
\begin{equation}
Q^{\pi}(\bstate_t, \mathbf{u}_t) = r(\bstate_t, \mathbf{u}_t) + \gamma Q^{\pi}(\bstate_{t+\Delta t},\mathbf{u}_{t+\Delta t}).
\label{eq:Bellman_quality}
\end{equation}

\medskip 
At this point, a few remarks are in order:
\begin{enumerate}
\item Since $\rho >0$ and $\Delta t>0$, the discount factor $\gamma \in (0,1)$. In the MDP framework, the value function can be represented as the cumulative discounted reward or \emph{return} $R_t$ defined by
\begin{equation}
\mathcal{J}^{\bpolicy}(\bstate_t) = R_t = \sum_{l =0}^\infty{\gamma^l r(\bstate_{t+l \, \Delta t})}.
\label{eq:sequence_problem}
\end{equation}
The two benefits of introducing a discounted reward is that the return is well defined for infinite series ($l \rightarrow \infty$), and that it gives a greater weight to earlier rewards, meaning that we care more about imminent rewards and less about rewards we will receive in the future.
\item Equation \eqref{eq:Bellman_opt} highlights how the discounted infinite-horizon optimal problem can be decomposed in a series of local optimal problems. This is the \textit{Bellman's principle of optimality}, \cite{bellman1957markovian}:  
An optimal policy has the property that whatever the initial state and initial decision are, the remaining decision must constitute an optimal policy with regard to the state resulting from the first decision. 
\item In \eqref{eq:Bellman_value}, the model \eqref{eq:statespaceNLN} does not appear explicitly. If $\bstate_t$ can be observed, and the reward $r$ can be measured, it is possible to recover $\mathcal{J}^\pi(\bstate_t)$ by inference from the agent-environment interaction ensuring that $\mathcal{J}^\pi(\bstate_t)$ is solution of the Bellman equation (\ref{eq:Bellman_value}). This is precisely what the \textit{Reinforcement Learning} (RL) approach does. 
In this framework, the policy and the value function are deterministic or stochastic functions.
These functions can be represented as tables when the number of states and actions are sufficiently low.
When there are more state and action variables, the \emph{curse of dimensionality} occurs \cite{bellman1961} and there is a need to approximate these functions as parameterized functional forms. Then, the \textit{learning} process consists in optimizing the parameters so that the value function is maximized and the constraint imposed by the Bellman optimality \eqref{eq:Bellman_opt} is satisfied.
In \emph{Deep Reinforcement Learning} (DRL), the functions are represented by a Neural Network (NN). 
\end{enumerate}

Within the RL framework, the Bellman equation is formulated in terms of the expectation of the value function. 
In the value Bellman equation \eqref{eq:Bellman_value}, the value function is replaced by the expected value ($\mathbb{E}\left[\mathcal{J}^{\bpolicy}(\bstate_{t})\right]$). Similarly, in the cost-value Bellman equation \eqref{eq:Bellman_quality}, the Q-function is replaced by $\mathbb{E}\left[Q^{\pi}(\bstate_{t},\mathbf{u}_t)\right]$.

% DDPG
\section{Reinforcement Learning} \label{sec:RL}

As mentioned above, the aim of RL is to find the optimal policy $\bpolicy^\star$ that maximizes the return \eqref{eq:sequence_problem}. A crucial aspect for our application is the possibility for RL algorithms of learning directly from the observations resulting from the interactions of the agent with the environment.
In that sense, we do not rely on the full knowledge of the state, but on a few localized measurements. Hereafter, the value function, the policy and the whole formulation will be given as parameterized functions of the observable $\bsens_t$. In Sec.~\ref{sec:plant}, pointwise measurements will be used for the application on the KS system. In the following, we introduce a possible classification of the RL algorithms and specify our choices.

\subsection{Several classes of RL algorithms} \label{Sec_severalRLs}
In the literature, many RL algorithms exist.
In practice, the choice of the right algorithm often depends on the type of available actuators and sensors, on the system to control and on the task to fulfill. 
A general classification that embeds all RL algorithms can be made. 
Three classes of algorithms can be identified: i) \emph{Actor--only}, ii) \emph{Critic--only} and iii) \emph{Actor--Critic}, where the words \emph{actor} and \emph{critic} are synonyms for the policy and value function, respectively. 

\medskip
Actor--only methods work with a parameterized family of policies. The gradient of the value function, with respect to the parameters, is estimated, and the parameters are updated in a direction of improvement. In the DRL approach, the policy is represented by a neural network to be inferred such that $\mathbf{u} = \bpolicy(\bsens\vert \bomega)$, with $\bomega$ the parameters of the NN. The learning process is guaranteed by choosing $\bomega$ such that the discounted reward $R_t$ is maximized. In such a procedure, a single policy is evaluated by recording the system for a long time.
This evaluation allows the computation of $R_t$, and of the gradient of the value function with respect to the NN parameters of the policy, \textit{i.e.}, $\nabla_{\bomega} \mathcal{J}^{\bpolicy}(\bsens_t)$.
The algorithms belonging to this class are referred to as \texttt{REINFORCE} algorithms, \cite{sutton1998introduction}. A gradient-based optimization of the policy representation is carried-out with Stochastic Gradient Descent (SGD) algorithms (or closely related SGD-like strategies, for a recent review see \cite{ruder2016overview}). A possible drawback of such methods is that the gradient estimations may have a large variance.

\medskip
Critic--only methods rely exclusively on value function approximation and aim at learning an approximate solution to the Bellman equation, which will then hopefully prescribe a near-optimal policy.
In the "Deep" flavor of the algorithms, the state-action value function $Q^{\bpolicy}$ is approximated by a NN, leading to $Q^{\bpolicy}\left(\bsens_t,\mathbf{u}_t\vert \mathbf{\btheta}\right)$, with $\mathbf{\btheta}$ being the parameter vector of the NN.   
The parameters $\btheta$ are chosen such that the function $Q^{\bpolicy}$ is solution of the Bellman equation \eqref{eq:Bellman_quality}. 
A Q-learning algorithm \cite{watkins1992q} may be used to approximate the optimal action-value function.
The core of the algorithm is to update iteratively the function by using a weighted average of the old value and the new information:
\begin{eqnarray}
Q^{\bpolicy}\left(\bsens_t,\mathbf{u}_t\vert \mathbf{\btheta}\right) 
&\longleftarrow&
Q^{\bpolicy}\left(\bsens_t,\mathbf{u}_t\vert \mathbf{\btheta}\right) \\
&+& \alpha
\left(
r(\bsens_t,\mathbf{u}_t)
+
\gamma\,
\max_{\mathbf{u}}Q^{\bpolicy}\left(\bsens_{t+\Delta t}, \mathbf{u}\vert \mathbf{\btheta}\right)
-
Q^{\bpolicy}\left(\bsens_t,\mathbf{u}_t\vert \mathbf{\btheta}\right)\right), \nonumber
\end{eqnarray}
where $\alpha$ is the learning rate ($0<\alpha\le 1$).

The use of NN in the Q-learning procedure establishes the Deep Q-Networks (DQN), \cite{mnih2015human}.
This class of algorithms is often referred in the specialized literature as the first \textit{"human-level"} control algorithm. The optimal policy does not need a function representation as it shall consist of choosing the action $\mathbf{u}$ that maximizes the optimal state-action value. For a finite number of possible actions, the Q-function is computed for each of these and the action that ensures the maximum is then chosen. The use of NN imposes the implementation of some strategies to regularize the learning process, such as memory, prioritized experience, target neural network, etc., \cite{mnih2015human}. 

\medskip
From a mathematical point of view, the difference between the Actor--only and Critic--only approaches relies on the difference between the Pontryagin's maximum principle and the Bellman principle. The first principle is based on a \textit{perturbative} approach: the neighborhood of the controlled trajectory is explored by perturbing the states in order to further minimize the cost function (curve optimization). The second principle is based on the solution of the HJB equation under the assumption of Markovianity of the system (function optimization) as explained in the previous paragraphs. While the Pontryagin's maximum principle is a necessary condition for the optimality, the fulfillment of the HJB equation is a necessary and sufficient condition if and only if the action-state is fully known (Legendre-Clebsch condition, see \cite{Choset-2005-9207}).

\medskip
The last class of algorithms, the Actor--Critic ones, combines the advantages of the two aforementioned approaches.
The critic uses an approximation architecture to learn a value function, which is then used to update the actor's policy parameters in a direction of performance improvement. 
In DRL, two NNs are employed simultaneously, one for the policy and a second one for the value function. In that sense, a possible implementation consists in combining the \texttt{REINFORCE} and the Q-learning procedures such that the optimal controller is obtained, \cite{mnih2016asynchronous}. From the algorithmic viewpoint, the coupling of the two NNs is due to the update of the networks that is performed simultaneously by sharing the same expectation of the discounted reward. With a parameterized policy function, the limitation of the application of the Critic--only approach to a discrete action space is circumvented.

\subsection{Deep Deterministic Policy Gradient as an actor-critic algorithm for RL} \label{sec:ddpg}
\begin{figure}[t]
\begin{center}
    \includegraphics[width=9.5cm]{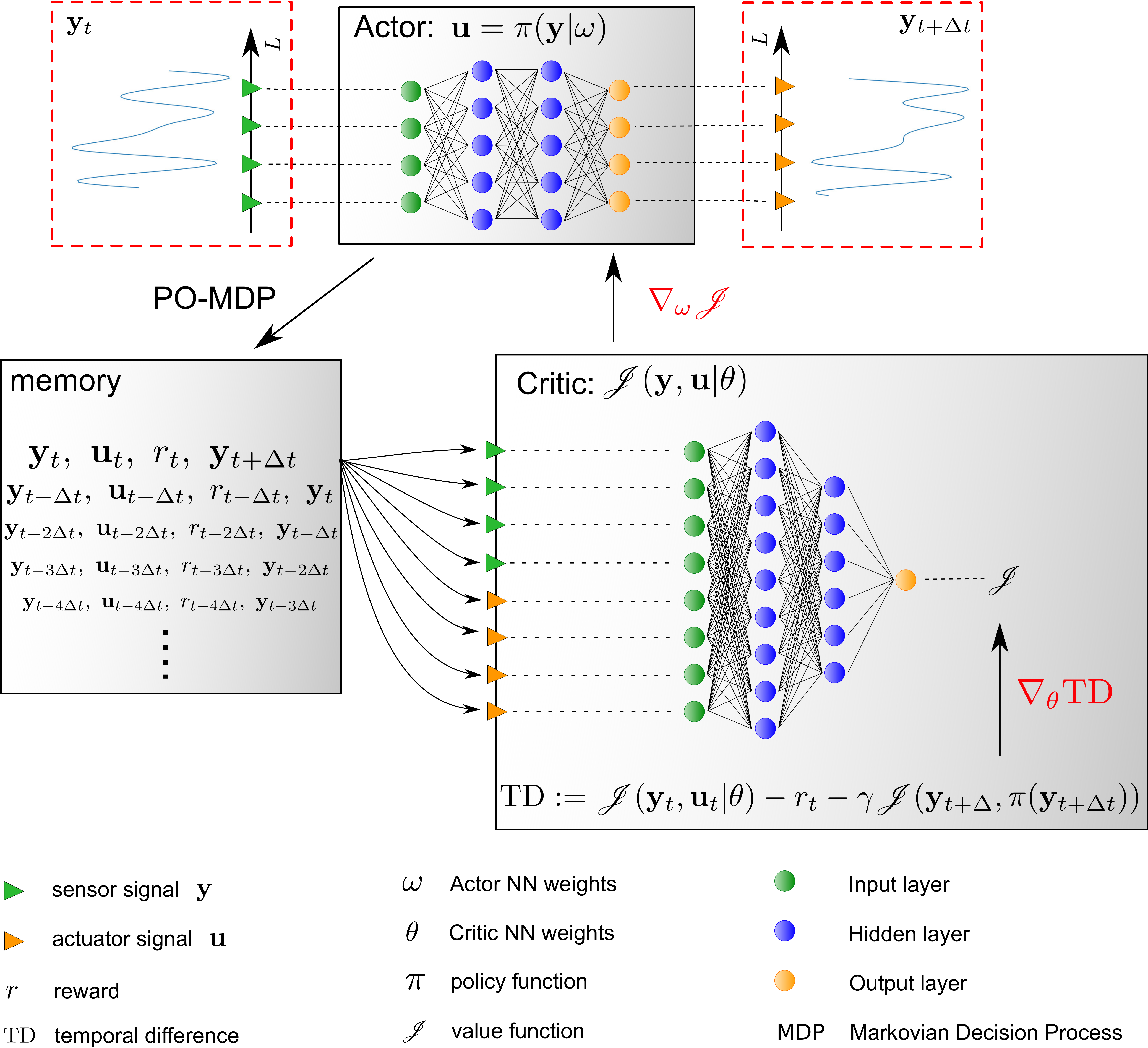}
    \caption{Graphical sketch for the Deep Deterministic Policy Gradient algorithm. The strategy is an Actor--Critic controller, where each of the part is implemented by means of a neural network, as depicted in the sketch.}
    \label{fig:DDPG_shetch}
\end{center}
\end{figure}
In this work, we consider the Deep Deterministic Policy Gradient (DDPG) algorithm \cite{lillicrap2015continuous} as an Actor--Critic, model free algorithm (see the graphical summary in Fig.~\ref{fig:DDPG_shetch}).
The advantage of DDPG compared to the Deep Q Network considered in \cite{mnih2015human} is to handle continuous action domain that are often encountered in physical control tasks.
Essentially, the DDPG algorithm
adapts the theoretically-grounded Deterministic Policy Gradient (DPG) algorithm introduced in \cite{silver2014deterministic} to the Deep Reinforcement Learning setting where NNs are employed to represent the policy and value functions. 
The basic idea of DPG is to update the policy parameter $\bomega$ in the direction of the gradient of $Q^{\bpolicy}$, rather than globally maximizing $Q^{\bpolicy}$ as it is done classically in policy gradient methods.
After application of the chain rule to $Q^{\bpolicy}\left(\bsens_t,\bpolicy(\bsens_t \vert\bomega)\right)$, we can show that the policy improvement at each iteration number $k$ can be decomposed into the gradient of the state-action value with respect to actions, and the gradient of the policy with respect to the policy parameters (see Eqs. 6 and 7 in page 3 of \cite{silver2014deterministic}; \textit{ibid.} the policy $\bpolicy$ is indicated as $\mu$). We finally obtain:
\begin{equation}
\bomega^{k+1}=
\bomega^k + 
\alpha 
\mathbb{E}
\left[
\bnabla_{\bomega}
\bpolicy(\bsens_t \vert\bomega)
\left.
\bnabla_{\mathbf{u}_t}
Q^{\bpolicy^k}\left(\bsens_t,\mathbf{u}_t\right)
\right|_{\mathbf{u}_t=\bpolicy(\bsens_t \vert\bomega)}
\right],
\label{eq:Silver_DPG}
\end{equation} 
where $\alpha$ is a learning rate.

The tuple corresponding to one time discrete transition $\{ \bsens_t, \mathbf{u}_t, r_t, \bsens_{t+\Delta t}\}$ defines a Markov Decision Process (MDP) or Partially Observable Markov Decision Process (PO-MDP) in case only a partial measurement of the state is accessible. At each iteration, the PO-MDP is stacked in memory and successively used by the critic optimizer to reduce the Temporal Difference error defined as
\begin{equation}
Q^{\bpolicy}\left(\bsens_t,\mathbf{u}_t\vert \mathbf{\btheta}\right)
-
r(\bsens_t,\mathbf{u}_t)
-
\gamma\,
Q^{\bpolicy}\left(\bsens_{t+\Delta t},\mathbf{u}_{t+\Delta t}\vert \mathbf{\btheta}\right).
\end{equation}
Regardless of the optimized critic NN, the actor NN is also optimized according to \eqref{eq:Silver_DPG}.

The optimality of the control is guaranteed by the Bellman principle under the hypothesis that the state-action space is known. For this reason, the state-action space needs to be explored. The exploration is carried out by perturbing the parameters of the policy as
\begin{equation}
\mathbf{u}_t = \bpolicy(\bsens_t \vert \bomega +\mathcal{N}_1) +\mathcal{N}_2,
\end{equation} 
where $\mathcal{N}_1$ and $\mathcal{N}_2$ are two noise processes. $\mathcal{N}_1$ is essential at the beginning of the exploration process to introduce unbiased and uncorrelated random actions. The NN is a strongly non-linear representation of the true optimal policy function and a small variation of the parameters can thus lead to drastic changes in the output action. Both noise processes $\mathcal{N}_1$  and $\mathcal{N}_2$ vary over time.  $\mathcal{N}_1$ is damped to a desired standard deviation of the resulting action $\mathbf{u}$, while the amplitude of $\mathcal{N}_2$ is given as a function of an Ornstein-Uhlenbeck process. This different choice introduces two different time-scales in the exploration process: when $\mathcal{N}_1$ is already damped, the noise process $\mathcal{N}_2$ still allows to explore efficiently the control landscape in the vicinity of the converged policy. We refer to DDPG articles for further technical information about its implementation, \cite{silver2014deterministic, lillicrap2015continuous, schaul2015prioritized}.

% KS - PhaseSpace - Results
\section{Control of the Kuramoto-Sivashinsky model}\label{sec:plant}
The DDPG algorithm introduced in the previous section is tested on the one-dimensional (1D) Kuramoto-Sivashinsky (KS) equation. The 1D KS equation is used for the description of flame fronts and reaction-diffusion systems. It is well known and referenced as being one of the simplest nonlinear PDEs exhibiting spatio-temporal chaos \cite{holmes1996turbulence,cvitanovic2010state}, thus providing an appropriate test-bed for the proposed control strategy. In the following, we first describe the dynamics of the KS in the phase-space (Sec.~\ref{sec:dynaKS}), before a description of the controlled cases is given in Sec.~\ref{sec:results}.
\begin{figure}
\centering
\includegraphics[width=11cm]{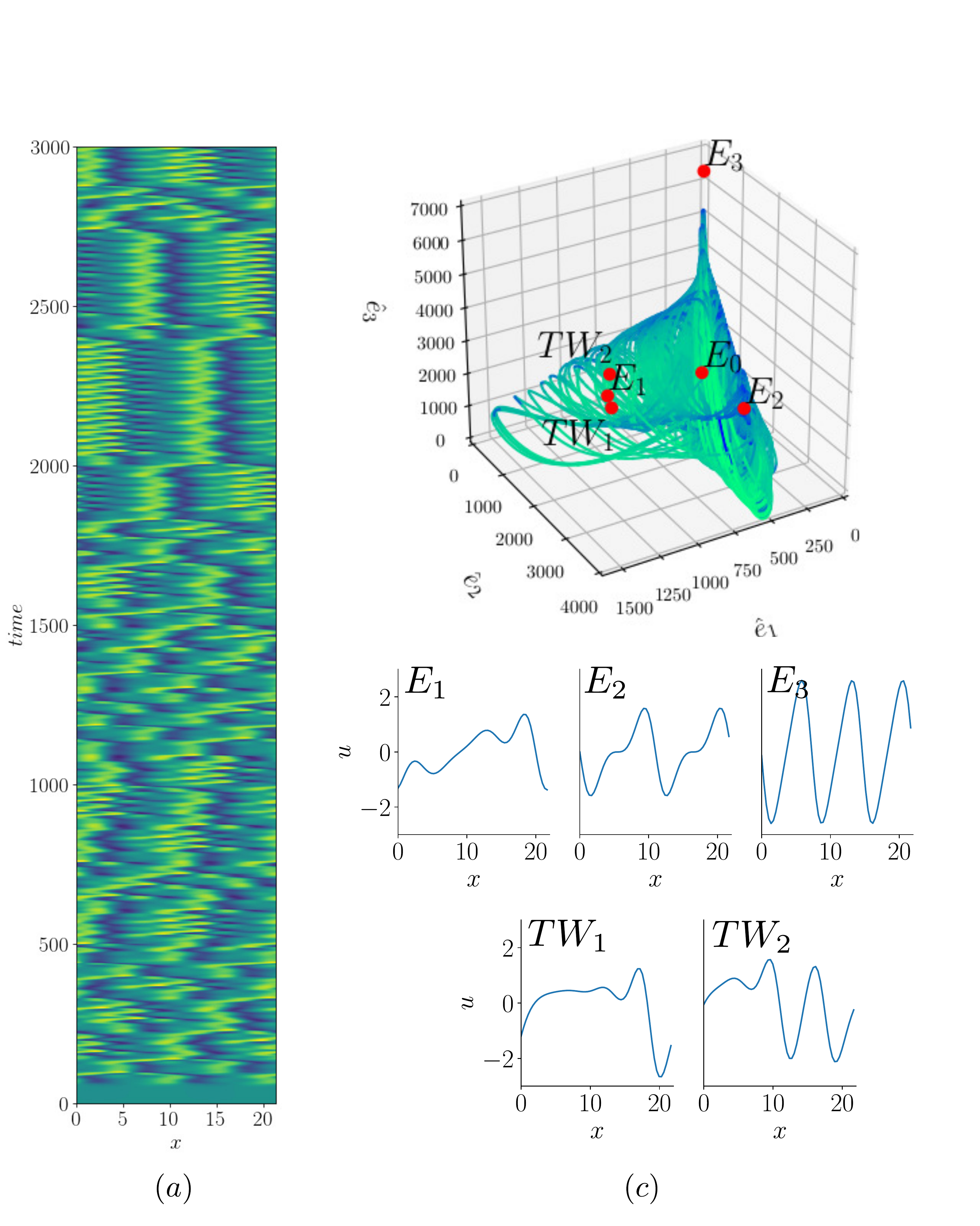}
\caption{%
Dynamics of the Kuramoto-Sivashinsky (KS) equation on a time interval of $3000$ time-units.
The contour plot in (a) shows the space-time behavior of a trajectory emanating from the perturbed state $E_0$.
The same trajectory is shown in (b) in the phase-space spanned by the dominant Fourier coefficients $\left\{\hat{e}_1,\hat{e}_2,\hat{e}_3\right\}$.
Red spots indicate: the trivial solution $E_0$, the three non-trivial invariant solutions $E_i$, $i=1, 2, 3$, and the two traveling waves $TW_1$ and $TW_2$. The spatial distributions of the solutions are given in (c).
Note how the trajectory is attracted by $E_3$, although never visited, while the other two equilibria are often visited.}
\label{fig:sim}
\end{figure}

\subsection{Dynamics of the Kuramoto-Sivashinsky equation}\label{sec:dynaKS}
The time evolution of the velocity $v=v(x,t)$ on a periodic domain of length $L$ is given by
\begin{equation}
\dfrac{\partial v }{\partial t}+v\dfrac{\partial v }{\partial x} 
= -\dfrac{\partial^2 v }{\partial x^2} -\dfrac{\partial^4 v }{\partial x^4} + g,
\label{eq:kuramoto}
\end{equation}
where $g$ is a spatio-temporal forcing term. The equation is characterized on the left-hand side by a nonlinear convective term and on the right-hand side by a diffusion term expressed by two addends: a $2nd$-order derivative related to energy production, and a $4th$-order derivative acting as an hyper-diffusion.

The length of the domain dictates different regimes for the solution. Introducing the trivial solution $E_0=0$, it can be shown that for $L<L_c = 2\pi$, the dynamics is stable and converges towards $E_0$, while for $L>L_c$ a chaotic dynamics emerges (in a Lyapunov sense). In this work, we focus on the dynamics of the KS equation for a domain length of $L=22$, identical to that studied in  \cite{cvitanovic2010state}. This length is large enough for exhibiting many of the features typical 
of turbulent dynamics observed in large KS systems, see \cite{pathak2018model}, but sufficiently small for a thorough analysis of the state-space dynamics as performed by \cite{greene1988steady,cvitanovic2010state}. 

\subsubsection{Numerical Simulations}
Numerical simulations are used for solving \eqref{eq:kuramoto}. 
As already mentioned, periodic boundary conditions are imposed, $v(x,t)= v(x+L,t)$. The spatial periodicity allows us to project the instantaneous solution onto Fourier modes and to perform spatial derivatives in the Fourier space. For $L=22$, $N=64$ Fourier collocation points have been used. 
This number of collocation points is deemed sufficient for an accurate representation of the spatio-temporal dynamics (see the related discussion in \cite{cvitanovic2010state}). Time marching is carried out with a 3rd-order semi-implicit Runge-Kutta scheme \cite{kar2006semi}. The linear operator on the right-hand side of \eqref{eq:kuramoto} is marched in time by an implicit scheme, whereas the non-linear and forcing terms are marched in time explicitly. For all numerical simulations, a time step of $0.05$ was adopted.

In Fig.~\ref{fig:sim}(a), a sample of the chaotic dynamics is shown
when the solution is initialized with the trivial state $E_0$ perturbed by a Gaussian white noise with an amplitude of order $10^{-4}$.
The trajectory evolves in time as shown in the corresponding phase-space \ref{fig:sim}(b) obtained by projecting the dynamics onto the set of the dominant Fourier modes $\hat{e}_1$, $\hat{e}_2$, $\hat{e}_3$. Indeed, due to the periodicity of the domain, the velocity fields are characterized by different phases, with respect of which the Fourier projection is independent. In  Fig.~\ref{fig:sim}(b), the different invariant solutions computed by Newton iterations are shown, namely the three unstable, fixed points indicated by $E_i$ ($i=1,2,3$) and the two traveling waves indicated by $TW_i$ ($i=1,2$). 
Apart for the relevant role in the dynamics, the invariant solutions $E_i$ are used as target of our control strategy in Sec.~\ref{sec:results}. A detailed analysis of the resulting dynamics shown in Fig.~\ref{fig:sim}(a)$-$(b) is beyond the scope of the present paper and already well described in \cite{cvitanovic2010state}. However, it is interesting to observe how the system evolves by frequently visiting the vicinity of the solutions $E_1$ and $E_2$, as well as the two traveling waves; on the other hand, the dynamics is ``attracted'' by the unstable solution $E_3$, along the manifold associated with it, without ever reaching it.

\begin{figure}
   \begin{minipage}[b]{0.5\textwidth}
    \includegraphics[width=7.25cm]{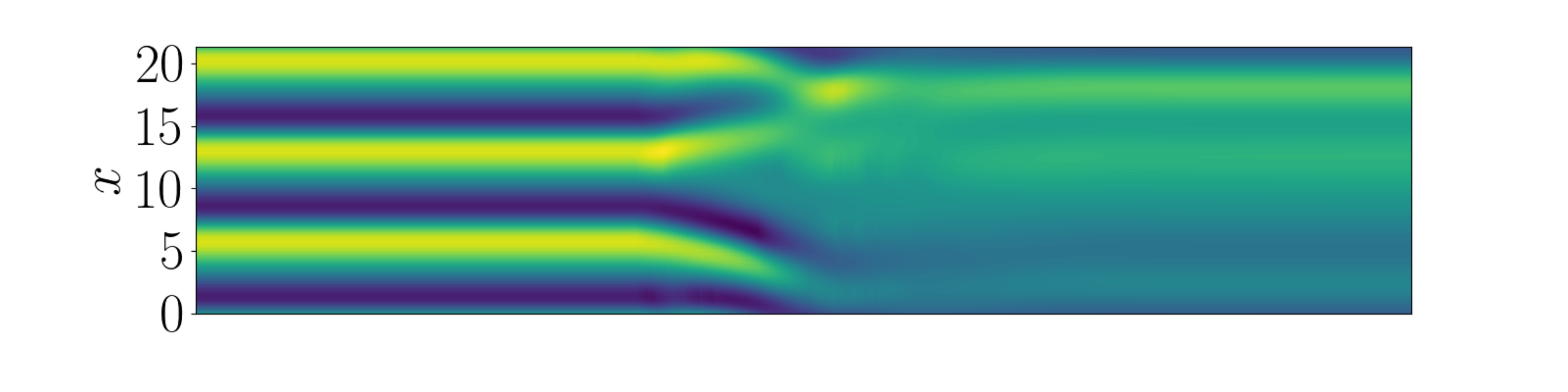} 
    \includegraphics[width=7.25cm]{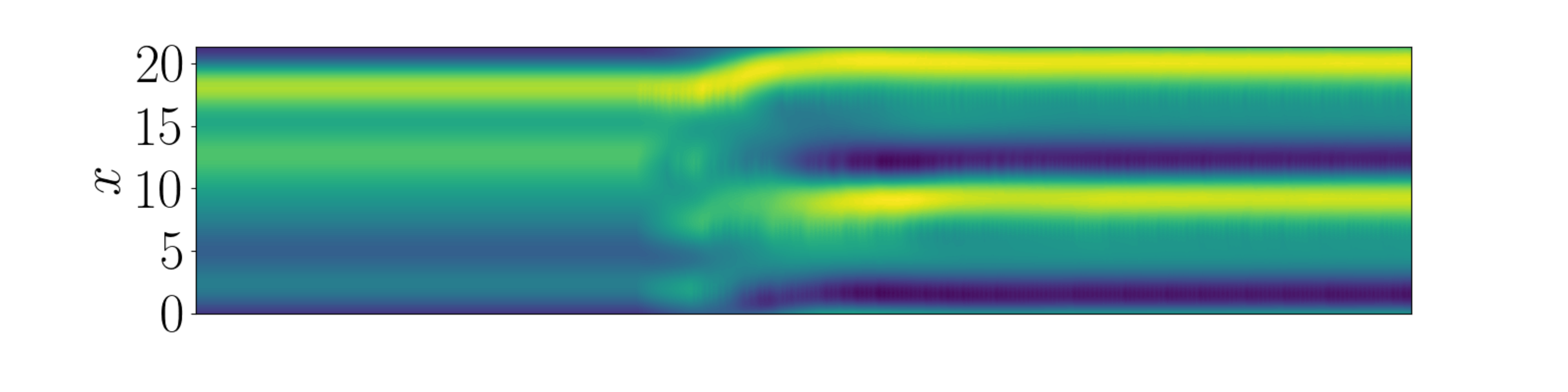} 
    \includegraphics[width=7.25cm]{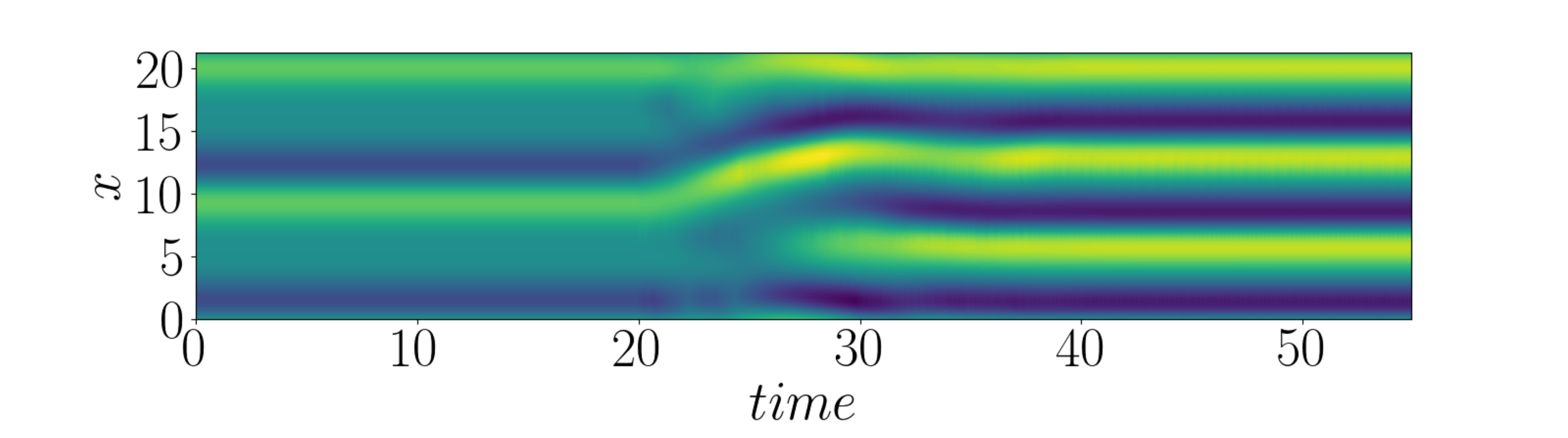} 
  \end{minipage}
  \begin{minipage}[b]{0.4925\textwidth}
    \includegraphics[width=6.75cm]{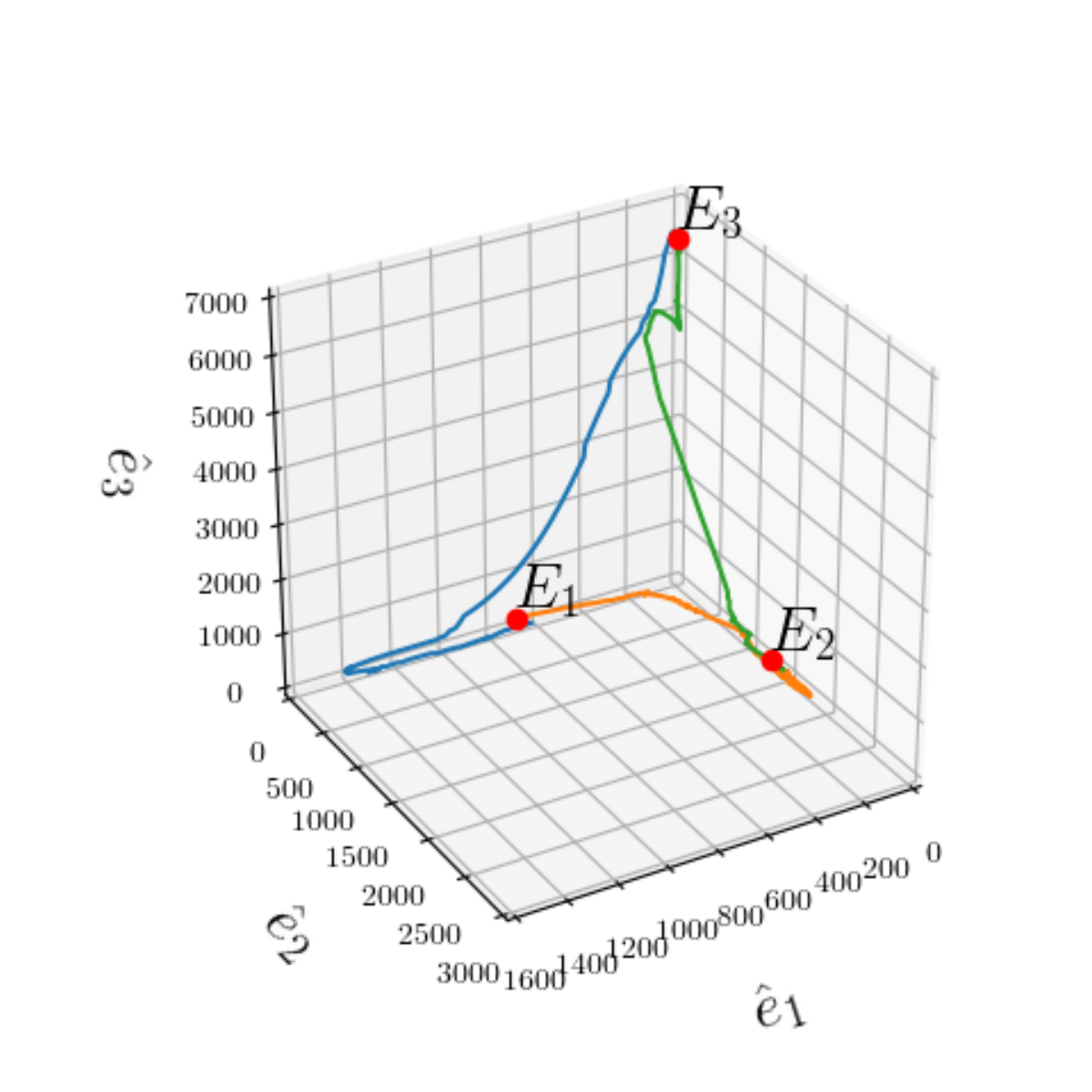} 
  \end{minipage}
        \put(-150,155){(d)}
        \put(-345,155){(a)}
        \put(-345,102.5){(b)}
        \put(-345, 50){(c)}
      \caption{Closed-loop dynamics for the three control test cases: (a) $E_3\rightarrow E_1$, (b) $E_1\rightarrow E_2$ and (c) $E_2\rightarrow E_3$. The controller is switched on at $t=20$. The control strategy used in each test case is the optimal RL policy. The trajectories are shown in the Fourier phase-space in (d).}
    \label{fig:con}
\end{figure}
\begin{figure}
\begin{center}
    \includegraphics[width=9.5cm]{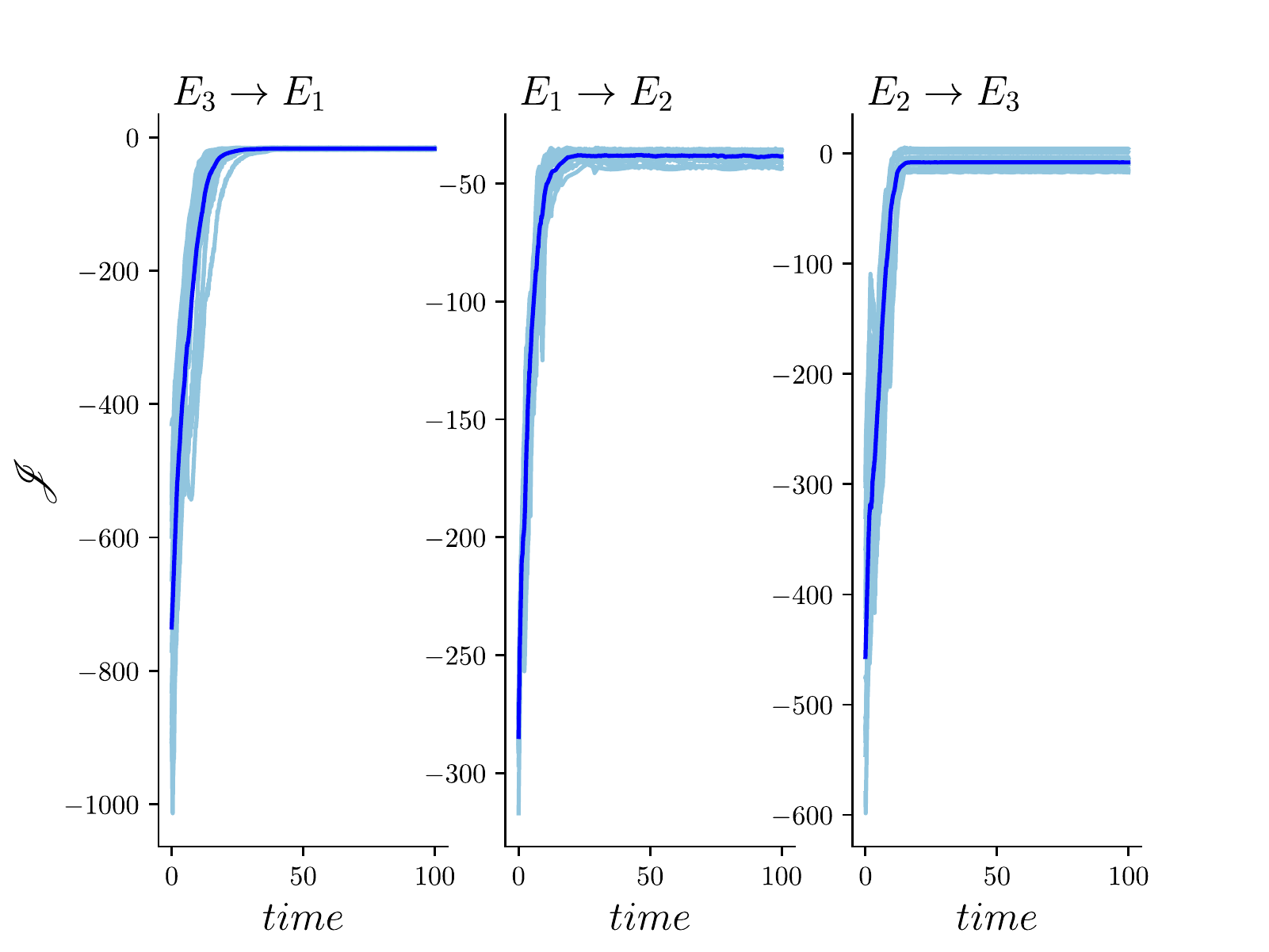}
      \caption{The value function is plotted for each of the three control test cases considered in Fig.~\ref{fig:con}. For each of the test case, the robustness of the policy is assessed by running $15$ different experiments emanating from different initial conditions of the control. The darker blue line indicates the average of the value as a function of time. For each case, it can be observed that - regardless of the initial conditions - the controller is capable of reaching the target in a robust manner.}
    \label{fig:value}
\end{center}
\end{figure}

\subsection{Controlled system}\label{sec:results}
As mentioned in Sec.~\ref{sec:nonopt}, the plant is the system to be controlled, including the inputs (actuators) and the outputs (sensors) of the system. The number of actuators and sensors, their distribution and the combination of the two define a rather large parametric space. In this work, we consider $4$ actuators equispaced along $x$. The support of each actuator is assumed Gaussian-shaped and the forcing term in \eqref{eq:kuramoto} is given by
\begin{equation}
g\left(x;\left\{x_i^a\right\}_i\right)=
\sum_{i=1}^4{\left(2\pi\sigma\right)^{-1\slash 2} \, 
\exp\left(-\dfrac{\left(x -x_i^a\right)^2}{2 \sigma^2}\right)},
\end{equation}
where $x_i^a\in \left\{0, L/4, L/2, 3L/4\right\}$ is the location of the actuators along the $x$-axis and $\sigma=0.4$ is the variance defining their spatial distribution. Regarding the sensors, we make the realistic assumption that real-life controllers rely only on partial information. We then  introduce $4$ equispaced sensors measuring the local velocity $v$. These sensors are staggered with respect of the actuator's locations and are located at $x_i^s\in\left\{L/8, 3L/8, 5L/8, 7L/8\right\}$.

\subsubsection{Implementation of the control policies}
The Deep Deterministic Policy Gradient  approach (see Sec.~\ref{sec:ddpg}) is implemented using scripts in \texttt{PyTorch}\footnote{https://pytorch.org/}. The critic and actor parts are both implemented by means of a neural network. 
The neural network of the critic part consists of an input layer of dimension $8$, with $4$ nodes dedicated to the actuator's signals and $4$ to the sensor's ones, and a single scalar output layer; two hidden layers are introduced, with $256$ and $128$ nodes, respectively, both featuring the \texttt{swish} activation function, \cite{ramachandran2017searching}. The neural network for the actor part approximates the control law. As such, the input layer is fed with the sensor measurements and is of dimension $4$, while the output layer provides the control signal feeding the $4$ actuators. Two hidden layers are used, of dimensions $128$ and $64$, with activation functions \texttt{swish} and \texttt{tanh}, respectively. The last layer acts as a saturating function. The maximum amplitude of the output is $g_\text{max}=0.5$.
The training is performed using the \texttt{Adam} optimization algorithm, with learning rate of $0.001$ for both the networks, and mini-batch with $200$ examples for the optimization of the critic-network. Finally, the discount factor is set to $\gamma=0.99$.

\subsubsection{Results}
We want to demonstrate the ability of the DDPG-based policies to drive the chaotic system towards the unstable fixed points and keep the dynamics in their vicinity, using the localized sensors and actuators discussed before. From the engineering viewpoint, this would correspond to leading a fluid system to a given operating condition only relying on some measurements of the environment.
Specifically, we introduce three control test cases, each of them targeting to minimize the distance between the state $v$ and the non-trivial invariant solutions, namely $E_1$, $E_2$ and $E_3$, respectively. 
In the following, for simplicity, we label the different policies with the invariant solution considered in their objective function. 

First, we consider each individual policy for driving the dynamics of the system from one invariant solution to another. In Fig.~\ref{fig:con}, the spatio-temporal behavior of the solution is shown in the inserts (a)-(c), and as a phase-space representation in (d). In each case, the RL controller is active for $t>20$. For each control case, the numerical experiments are repeated $15$ times with different initial conditions of the control, to assess the robustness of the policy.
In  Fig.~\ref{fig:value}, we report the value as a function of time for each experiment (light-blue). The darker line indicates the average of the runs. By comparing Fig.~\ref{fig:con}(d) and Fig.~\ref{fig:value}, we can observe that the controller is capable of driving the system in about 10 time-units towards the target solution. We also show that the control action of the policies is robust since the general behavior is comparable for all the runs considered.

\begin{figure}[t]
   \begin{minipage}[b]{0.45\textwidth}
    \includegraphics[width=7.25cm]{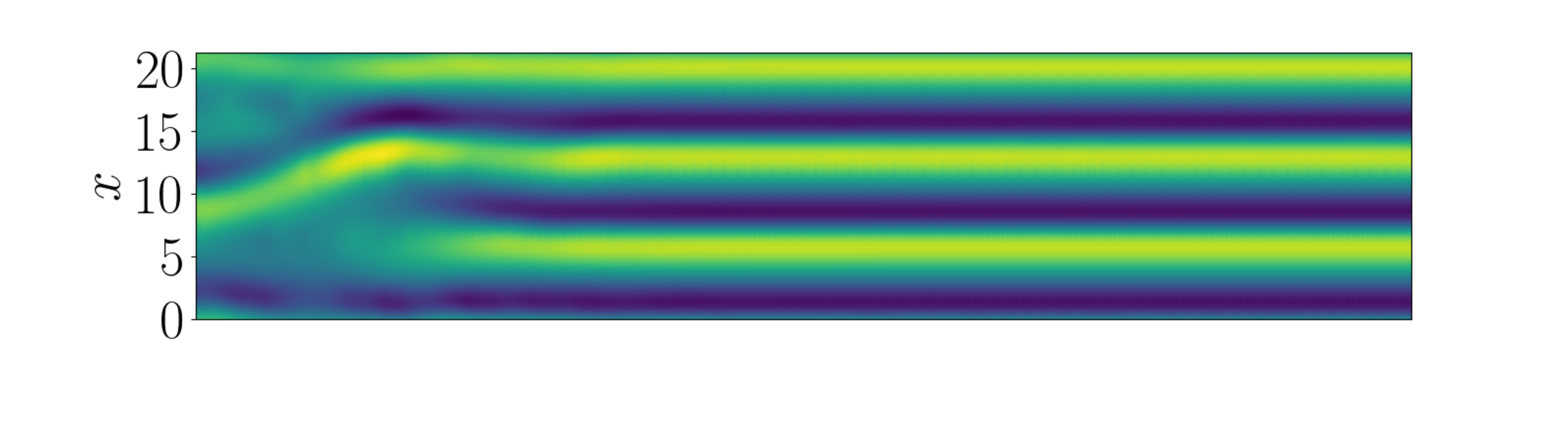} 
    \includegraphics[width=7.25cm]{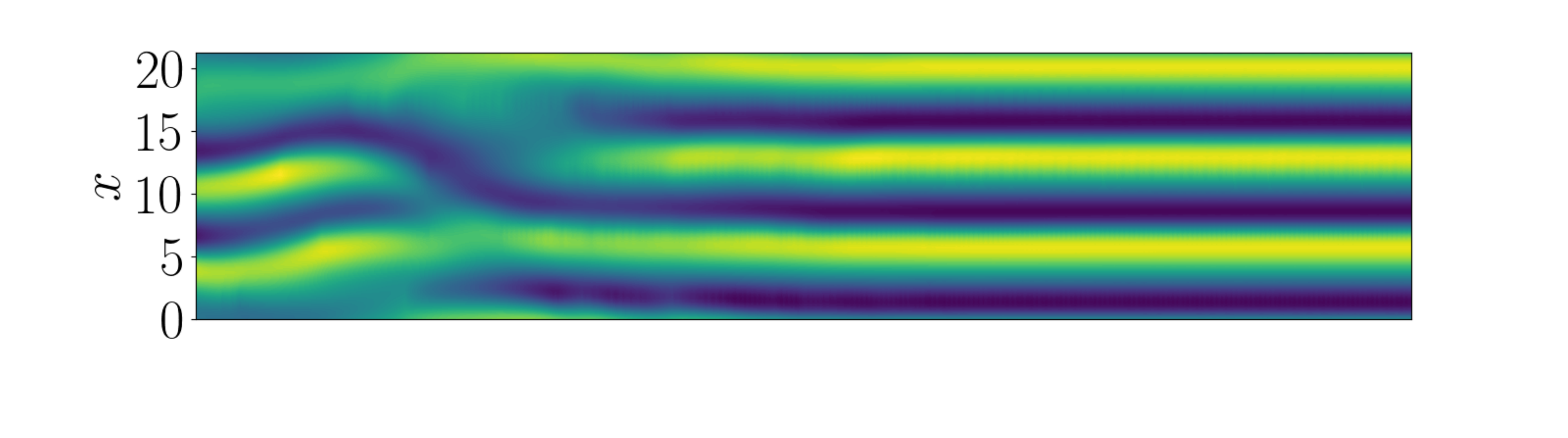} 
    \includegraphics[width=7.25cm]{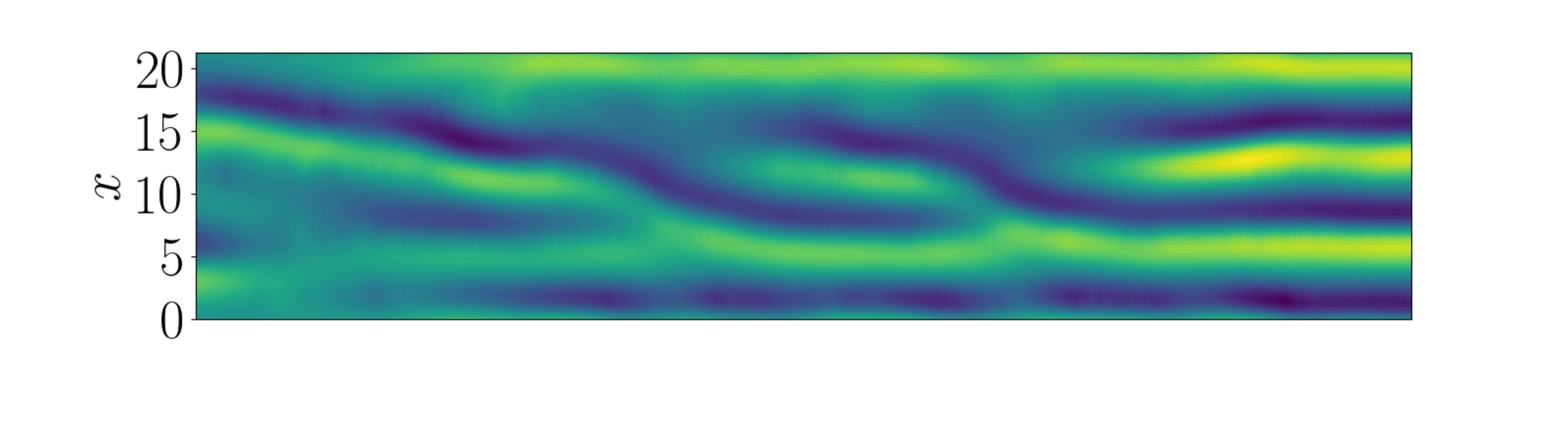} 
    \includegraphics[width=7.25cm]{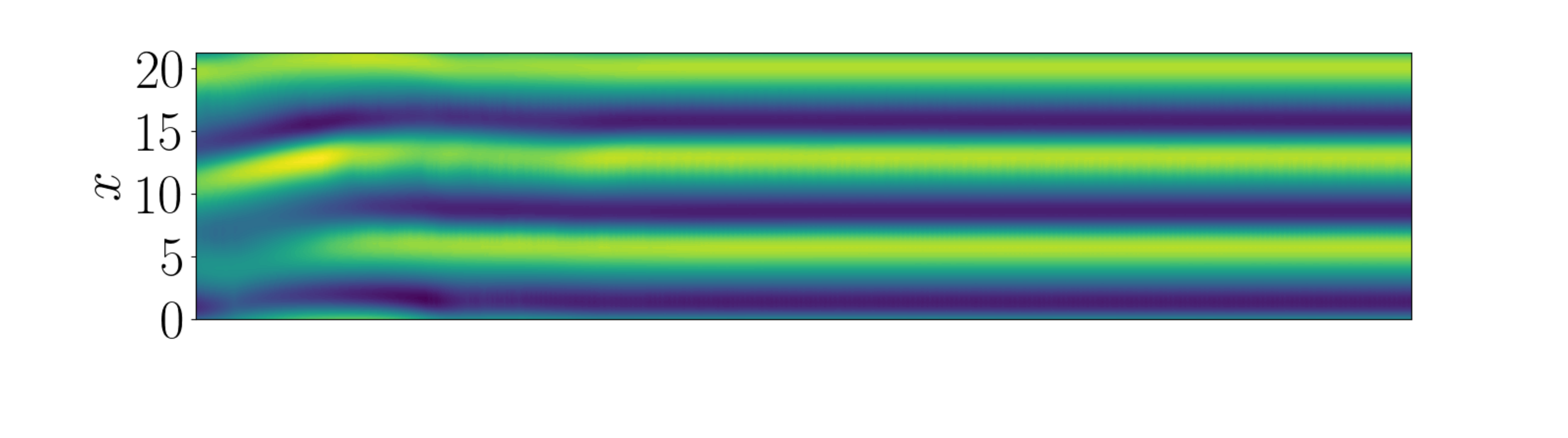} 
    \includegraphics[width=7.25cm]{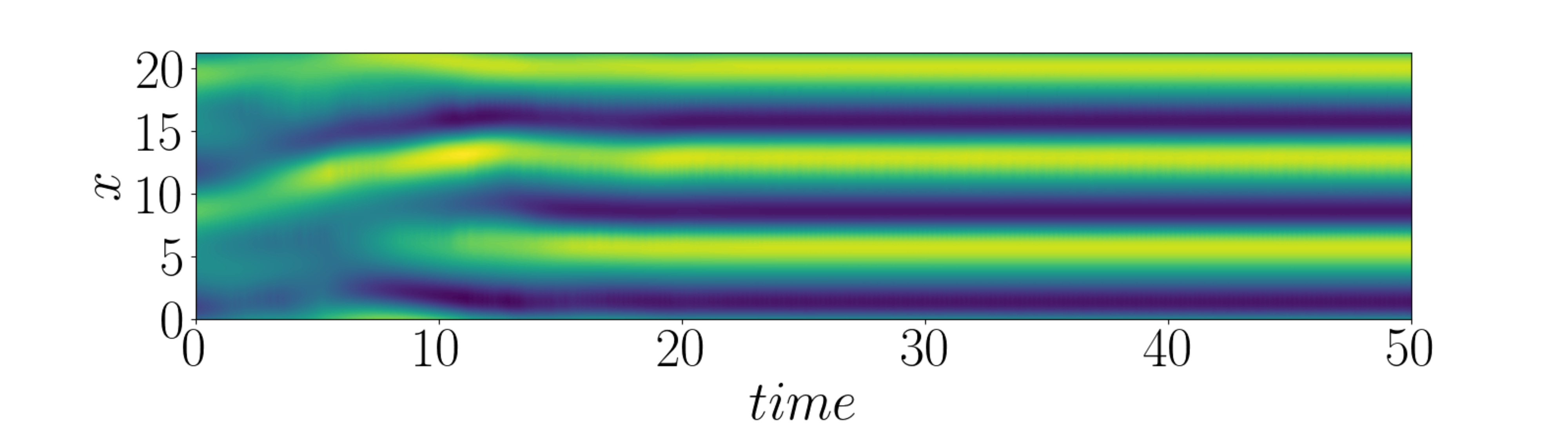} 
  \end{minipage}
  \begin{minipage}[b]{0.45\textwidth}
  \centering
    \includegraphics[width=6.75cm]{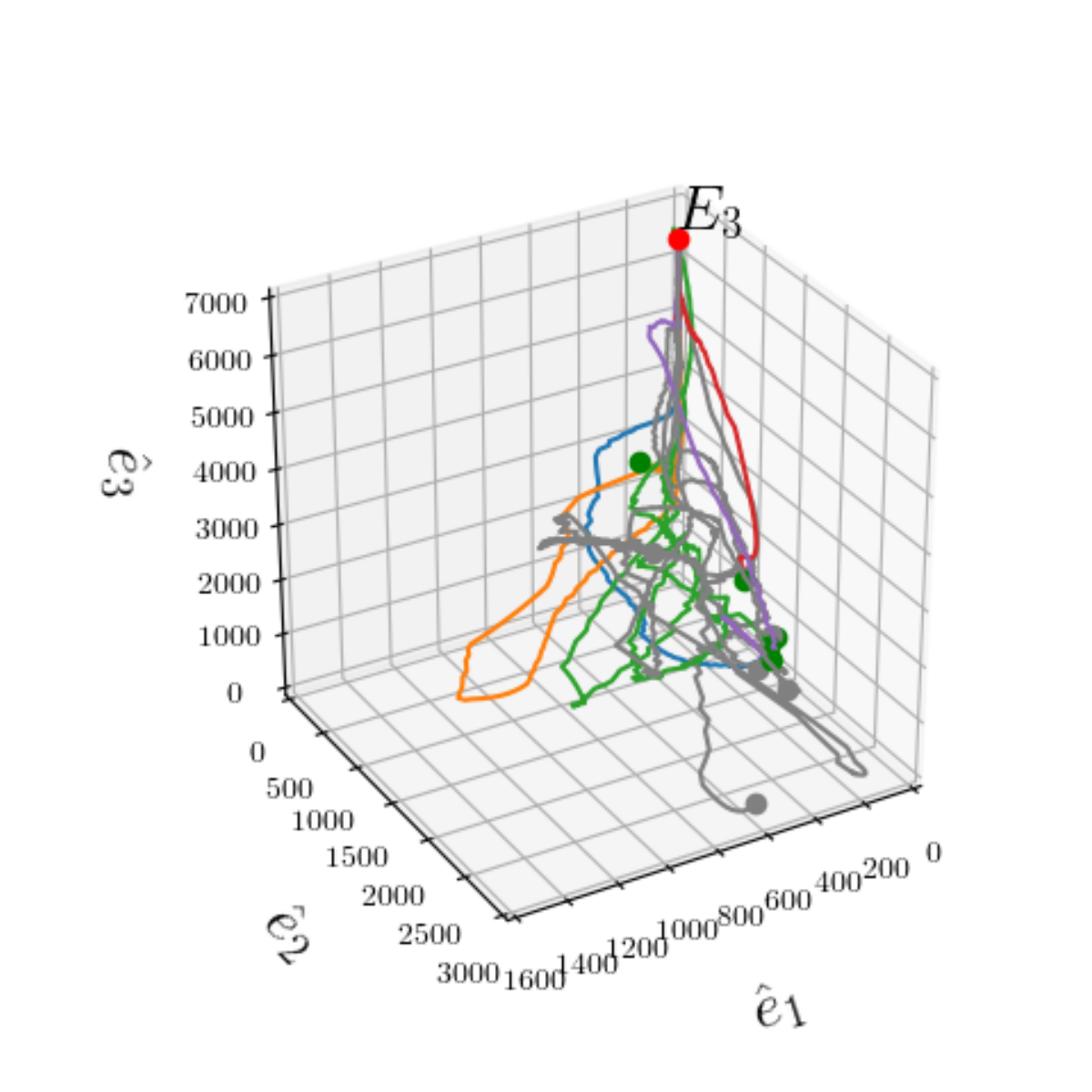} 
    \includegraphics[width=5.75cm]{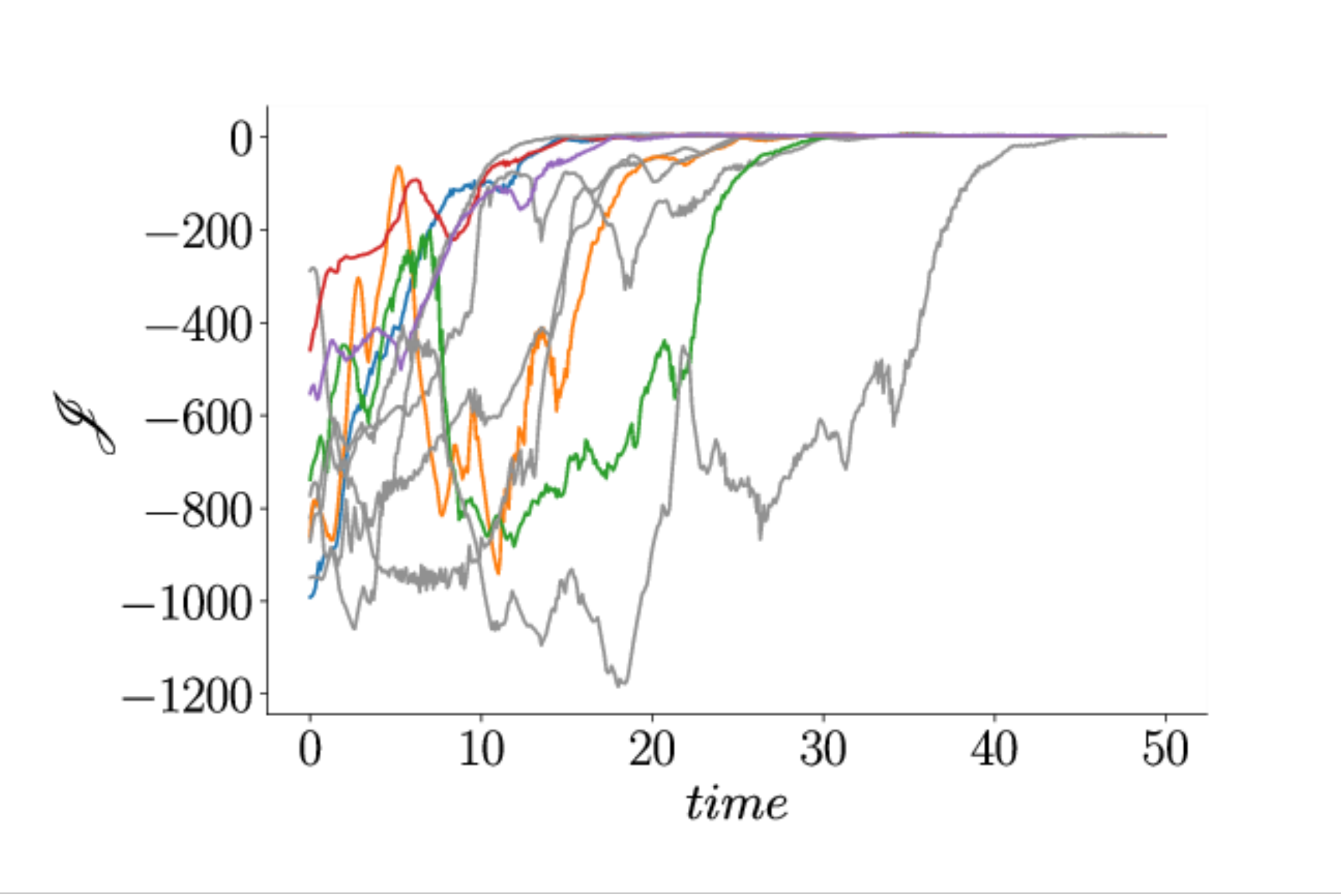} 
  \end{minipage}
        \put(-145,270){(f)}
        \put(-145,105){(g)}
        \put(-315,   50.0){(e)}
        \put(-315, 107.5){(d)}
        \put(-315, 167.5){(c)}
        \put(-315, 227.5){(b)}
        \put(-315, 287.5){(a)}
      \caption{Illustration of the robustness of the closed-loop dynamics with respect to changes of the initial conditions. Ten different initial conditions are randomly chosen for the optimal policy driving the dynamics towards the unstable state $E_3$. Five of these examples are reported in the figures (a)-(e). It is shown that, except for the trajectory shown in (c), the controller is capable to control the system in approximately $15$ time-units. The corresponding phase-space is shown in (f), where the trajectories in (a)-(e) are shown in colors, while the others are in gray. The corresponding values are shown in (g) as a function of time.}
    \label{fig:conrob}
\end{figure}

To further verify the robustness of the policies, we consider a second set of numerical experiments, where the initial condition is randomly chosen in the phase space. This is shown in Fig.~\ref{fig:conrob}, where we consider $10$ different initial conditions randomly and run the simulation under the control policy $\pi_{E_3}$. Five of these examples are reported in the inserts (a)-(e). The controller is capable of controlling the system in about $15$ time-units in all but one cases (c) where it took about 40 time-units. The phase-space view is shown in (f), where the trajectories in (a)-(e) are shown in colors, while the others are in gray. The time evolutions of the corresponding values are shown in (g). Large excursions of $\mathcal{J}$ roughly correspond to fluctuations in the phase-space, although all trajectories finally converge towards $E_3$. In that sense, the high independence from the initial conditions -- as compared to linear control methods such as the model-based LQG -- and the ability to maximize the value function demonstrate the level of performance and robustness of the DDPG-based strategy applied to the KS system.

% Conclusion
\section{Conclusions}\label{sec:conclusions}

We have presented a deep reinforcement learning (DRL) algorithm successfully controlling the chaotic dynamics governed by the well-known nonlinear Kuramoto-Sivashinsky (KS) equation. The DRL combines reinforcement learning principles with deep Neural Networks for approximating the value function and the control policy. Among the different available strategies, we have tested an actor-critic algorithm, the Deep Deterministic Policy Gradient (DDPG). This choice is related to the possibility of directly tailoring the classical strategies from the nonlinear optimal control theory and, in particular, the solution of the Hamilton-Jacobi-Bellman equation, with the resulting approximation obtained by DDPG. From this perspective, our approach departs from recent efforts found in the flow-control literature as it allows knowledge of the cost landscape in the vicinity of the system trajectory, hence bringing robustness with respect to perturbations, while allowing for high control performance. Further, our off-policy framework is not episode-based and can then potentially require less training data (faster learning).

\smallskip
We emphasize that, in the present work, the controlled system exhibits a chaotic behavior in the analyzed regime. This configuration is significantly more involved than systems exhibiting a periodic, or quasi-periodic, dynamics. Nonetheless, we demonstrate that, using a limited, localized set of actuators and sensors and model-free DRL controllers, it is possible to drive and stabilize the dynamics of the KS system around its unstable fixed solutions, here considered as target states. Also, we show that the controllers are robust with respect of the initial conditions by considering numerous trajectories emanating from them; for all the tested cases, the performances are essentially unchanged with respect of the initial conditions, in contrast with the application of more classic, linear controllers where this choice is often critical.

\smallskip
The complexity of the KS system and the possibility of computing a DRL control policy by solely using local measurements of the system suggest the extension of this application to the control of more realistic configurations such as turbulent boundary layers or mixers. In this sense, this work represents a first effort toward these applications. A number of developments are currently made to apply RL to these configurations addressing -- among other limitations -- limited and noisy observables, low measurement sampling frequency, high-dimensional state-action space and unknown and non-stationary time-delays.

%====================================================================
% Final things required by the template
%====================================================================

% \ethics{Not appliable.}

% \dataccess{No access given to GitLab}

\bigskip
\paragraph*{Authors' contributions}
{The project has been initiated by LM and LC. MAB implemented all the routines used in this work for the analysis of the dynamics governed by the KS equation, the approximation and control by RL, with supervision from OS and LM, and feedbacks from AA. The paper was written by MAB, OS, LM and LC, with feedbacks from GW and AA.}

\medskip
\paragraph*{Competing Interests}{We declare we have no competing interests linked to this study.}

\medskip
\paragraph*{Acknowledgements}{This work has benefited from discussions with Charles Pivot whom is gratefully acknowledged. The authors are also thankful to Sylvain Caillou for his support in the numerical implementation of the algorithm.}

\medskip
\paragraph*{Funding}{This project was generously funded by the French Agence Nationale pour la Recherche (ANR) and Direction G\'en\'erale de l'Armement (DGA) via the FlowCon project (ANR-17-ASTR-0022).}

% \disclaimer{Not appliable.}

%====================================================================

\appendix
\label{sec:appA}
\section{Solving the linear optimal problem: Variational and HJB approaches}
In this section we briefly summarize how a classic Linear Quadratic Regulator (LQR) problem can be derived from the more general formulation of the Hamilton-Jacobi-Bellman equation developed in Sec.~\ref{sec:HJB}. This Appendix is not meant to be exhaustive, but instead puts in perspective the application of the Reinforcement Learning with respect to standard approaches from control theory applied in the last decade in flow control. For more details on the optimal control theory, the interested reader can refer to \cite{lewis2012optimal} and to the reviews \cite{bewley2007linear,fabbiane2014adaptive,sipp2016amr}.

% Before proceeding, however, it is  of interest to report the linear form of the state-space model, that can be rewritten under the hypothesis of a time-invariant system as
To begin with, we consider the assumption of a linear, time-invariant, state-space model
\begin{subequations}
\begin{align}
\dfrac{\ddroit \bstate}{\ddroit t} &= {A}\, \bstate +{B} \, \mathbf{u},\label{eq:statespaceLN_1}\\
\bsens &= {C} \, \bstate +{D} \, \mathbf{u}.\label{eq:statespaceLN_2}
\end{align}
\label{eq:statespaceLN}
\end{subequations}
In this case, the system matrix $A\in\mathbb{R}^{\nstate \times \nstate}$ is obtained from the discretization in space of the analyzed model, including boundary conditions; the spatial distribution of the $m$ actuators is indicated by the matrix $B\in\mathbb{R}^{N\times m}$, while $C\in\mathbb{R}^{\nobs \times \nstate}$ is the matrix including $\nobs$ sensors and $D\in\mathbb{R}^{\nobs \times m}$ is the feedthrough (or feedforward) matrix. The aim of LQR is to define a control signal $\mathbf{u}$, such that the quadratic objective function
\begin{equation}
\mathcal{J}(\bstate_t,\mathbf{u}_t) = \dfrac{1}{2}
\int_0^T {\left(\bstate^\transpose {Q}_1\bstate + \mathbf{u}^\transpose {Q}_2\mathbf{u} \right) \, \ddroit \tau}
\label{eq:obj_LN}
\end{equation}
is minimized, \cite{lewis2012optimal}. 
The first term of \eqref{eq:obj_LN} is related to the energy of the state $\bstate$, while the second one penalizes the energy expense of the control action. In this definition, the matrices $Q_1 \succcurlyeq 0 \in \mathbb{R}^{\nstate \times \nstate}$ and $Q_2 \succ 0 \in \mathbb{R}^{m \times m}$ respectively represent state and control penalties.

The objective is then to determine $\mathcal{J}^\star=\min_{\mathbf{u}}\mathcal{J}$ under the constraints of the state equation \eqref{eq:statespaceLN_1}. 
This optimal control problem can be solved by minimizing the augmented Lagrangian
\begin{equation}
%\mathcal{\widetilde{J}} =
\mathcal{L} = 
\dfrac{1}{2}\int_0^{T} {\left( \bstate^\transpose {Q}_1\bstate +\mathbf{u}^\transpose {Q}_2\mathbf{u}\right) \, \ddroit \tau} + 
\int_0^{T} {\mathbf{p}^\transpose\left(\dot{\bstate}
-{A} \, \bstate -B \, \mathbf{u}
\right) \, \ddroit \tau},
\end{equation}
where the governing equations with initial conditions $\bstate=\bstate_0$ act as constraint and the \emph{co-state} or \emph{adjoint state} $\mathbf{p}$ can be interpreted as a Lagrangian multiplier. The \emph{optimality system} is then obtained by imposing the optimality conditions on $\mathcal{L}$. It leads to the following system of coupled equations:
\begin{subequations}
\begin{align}
\mathcal{L}_{\mathbf{p}}=\bzero \quad & \Longrightarrow &
\dfrac{\ddroit \bstate}{\ddroit t} &= A \, \bstate +B \, \mathbf{u}, 
\qquad\qquad \bstate(0)=\bstate_0,
&\quad\text{Direct equation}
\label{eq:direct}
\\
\mathcal{L}_{\bstate}=\bzero \quad & \Longrightarrow &
\dfrac{\ddroit \mathbf{p}}{\ddroit t} &=-A^\transpose \, \mathbf{p} + Q_1 \, \bstate,
\qquad \mathbf{p}(T)=\mathbf{0},
&\quad\text{Adjoint equation}
\label{eq:adjoint}
\\
\mathcal{L}_{\mathbf{u}}=\bzero \quad & \Longrightarrow &
\mathbf{0} & = B^\transpose\mathbf{p} +Q_2\mathbf{u}.
&\quad\text{Optimality condition}
\label{eq:optimality}
\end{align}\label{eq:optimalitysystem}
\end{subequations}
The direct-adjoint system \eqref{eq:direct}-\eqref{eq:adjoint} is solved iteratively, by marching forward in time the direct equation, and backward in time the adjoint equation, forced by the second term on the right-hand side of \eqref{eq:adjoint}. The optimal condition is updated using the gradient given by \eqref{eq:optimality}. It can be shown that the optimality system can be directly solved by means of an associated continuous time algebraic Riccati equation (CARE), \cite{lewis2012optimal}.

We have shown in Sec.~\ref{sec:RL} that the HJB equation determines the minimum value function for a given dynamical system with an associated reward. This equation is directly linked to the class of optimal control problem treated previously. For this reason, we aim at obtaining the CARE starting from the HJB equation rewritten for the linear system as
\begin{align}
- \dot{\mathcal{J}^\star}(\bstate,t) & = \min_\mathbf{u} \mathcal{H}\left(\bstate(t),\mathbf{u}(t), \mathcal{J}_{\bstate}^{\star},t\right) \\
& = 
\min_\mathbf{u} \left[
\dfrac{1}{2}
\left(\bstate^\transpose {Q}_1 \, \bstate
+\mathbf{u}^\transpose {Q}_2 \, \mathbf{u} \right)
+\mathcal{J}_{\bstate}^{\star\transpose}\, \left({A} \, \bstate +{B} \, \mathbf{u}\right)
\right],
\label{eq:HJB_lin}
\end{align}
where $\mathcal{H}$ by definition is the Hamiltonian. The right-hand side of \eqref{eq:HJB_lin} is minimized for $\mathcal{H}_{\mathbf{u}} = \bzero$, leading\footnote{Since ${Q}_2$ is a positive definite matrix, the existence of ${Q}_2^{-1}$ is guaranteed.
} to the optimal control 
\begin{equation}
\mathbf{u}^\star=-{Q}_2^{-1} {B}^\transpose\mathcal{J}^\star_{\bstate}.
\label{eq:optimal_control_HJB}
\end{equation}
Plugging this expression in~\eqref{eq:HJB_lin}, we get
\begin{equation}
- \dot{\mathcal{J}^\star} = 
\dfrac{1}{2}\bstate^\transpose {Q}_1 \, \bstate
-\dfrac{1}{2} \mathcal{J}_{\bstate}^{\star \transpose} B \, {Q}_2^{-1} \, B^\transpose \mathcal{J}_{\bstate}^\star
+\mathcal{J}_{\bstate}^{\star \transpose} {A} \bstate.
\label{eq:HJB_linmin}
\end{equation}
At this point, the value function $\mathcal{J}^\star$ can be expressed as a function of the unknown, real symmetric positive-definite matrix $S(t)$ 
\begin{equation}
\mathcal{J}^\star(\bstate(t),t)=\dfrac{1}{2}\bstate^\transpose(t) S(t)\, \bstate(t).
\end{equation}
This expression can be derived in time $\dot{\mathcal{J}^\star} = \dfrac{1}{2} \bstate^\transpose \dot{S}\,\bstate$, and with respect to the state $\bstate$, leading to the vector\footnote{Note also that if we compare the optimal control \eqref{eq:optimal_control_HJB} determined with the Hamiltonian and the one that can be deduced from the optimality condition  
\eqref{eq:optimality}, we arrive to $\mathbf{p}=\mathcal{J}^\star_{\bstate}=S \, \bstate$.
} $\mathcal{J}^\star_{\bstate} = S\bstate$. Introducing this term in \eqref{eq:HJB_linmin}, and considering that only the symmetric part of $S(t)A$ contributes to the result, we obtain
\begin{equation}
-\bstate^\transpose \dot{S}(t) \bstate
=\bstate^\transpose {Q}_1 \, \bstate
-\bstate^\transpose S(t) \, B \, {Q}_2^{-1} \, B^\transpose S(t) \, \bstate
+\bstate^\transpose \left( S(t) {A} +A^\transpose S(t) \right) \bstate,
\end{equation}
from which it is possible to obtain the differential Riccati equation since the expression above is valid $\forall \bstate$:
\begin{equation}
-\dot{S}(t)= S(t){A} +A^\transpose S(t) -S(t) \, B \, {Q}_2^{-1}\, B^\transpose S(t) +Q_1.
\end{equation}
The optimal action is finally obtained as 
\begin{equation}
\mathbf{u}^\star=-{Q}_2^{-1} {B}^\transpose S \, \bstate =: K \, \bstate
\end{equation}
where $K(t)$ is the optimal control gain obtained as solution to the LQR problem. The corresponding time-invariant optimal control gain is obtained for $T \rightarrow\infty$, leading to the continuous time algebraic Riccati equation.
%
%\begin{equation}
%0 = S A + A^\transpose S -S B Q_2 B^\transpose S +{Q}_1\label{eq:Riccati}
%\end{equation}
In perfect analogy, it is possible to define an optimal estimation problem that, using the LQR and thanks to the separation principle, leads to the Linear Quadratic Gaussian compensator \cite{glad2000control,lewis2012optimal}.

As a conclusion, some observations can be made. First of all, since the Riccati equation is not directly solvable for systems characterized by a large number of degrees of freedom, say $N>10^4$, the direct-adjoint system \eqref{eq:direct}-\eqref{eq:adjoint} is often solved instead, see for instance \cite{bewley2016methods} and \cite{luchini2014adjoint}. Optimization based on a finite sliding temporal window leads to the class of Model Predictive Controllers, often characterized by a nonlinear loop; examples in fluid mechanics are provided by the seminal work by\cite{bewley2001dns} and the recent application by \cite{cherubini2013nonlinear}.

%\section*{References}
\bibliographystyle{elsarticle-num}
\bibliography{rspa2019}

\end{document}